\documentclass[11pt]{article}
\usepackage{amsfonts}
\usepackage{bbm}
\usepackage{mathrsfs}
\usepackage{amssymb}

\def\hang{\hangindent\parindent}
\def\textindent#1{\indent\llap{#1\enspace}\ignorespaces}
\def\re{\par\hang\textindent}

\title{On (De)homogenized Gr\"obner Bases\thanks{Project supported by the National Natural
Science Foundation of China (10971044).}}

\author{Huishi Li\thanks{e-mail: huishipp@yahoo.com}\quad\quad Cang Su\\
{\small Department of Applied Mathematics}\\
{\small College of Information Science and Technology}\\
{\small Hainan University}\\
{\small  Haikou 570228, China}}

\date{}

\begin{document}
\maketitle
\begin{center}
\begin{minipage}{120mm}
{\small {\bf Abstract.} In a relatively extensive context,
(de)homogenized Gr\"obner bases are studied systematically. The
obtained results reveal further applications of Gr\"obner bases to
the structure theory of algebras.}
\end{minipage}\end{center}{\parindent=0pt\vskip 6pt
{\bf 2000 Mathematics Classification} Primary 16W70; Secondary 68W30
(16Z05).\vskip 6pt {\bf Key words} Graded algebra, graded ideal,
Gr\"obner basis, homogenization, dehomogenization }

\def\QED{\hfill{$\Box$}} \def\NZ{\mathbb{N}}
\def \r{\rightarrow}

\def\normalbaselines{\baselineskip 24pt\lineskip 4pt\lineskiplimit 4pt}
\def\mapdown#1{\llap{$\vcenter {\hbox {$\scriptstyle #1$}}$}
                                \Bigg\downarrow}
\def\mapdownr#1{\Bigg\downarrow\rlap{$\vcenter{\hbox
                                    {$\scriptstyle #1$}}$}}
\def\mapright#1#2{\smash{\mathop{\longrightarrow}\limits^{#1}_{#2}}}
\def\mapleft#1#2{\smash{\mathop{\longleftarrow}\limits^{#1}_{#2}}}
\def\mapup#1{\Bigg\uparrow\rlap{$\vcenter {\hbox  {$\scriptstyle #1$}}$}}
\def\mapupl#1{\llap{$\vcenter {\hbox {$\scriptstyle #1$}}$}
                                      \Bigg\uparrow}
\def\v5{\vskip .5truecm}
\def\T#1{\widetilde #1}
\def\OV#1{\overline {#1}}
\def\hang{\hangindent\parindent}
\def\textindent#1{\indent\llap{#1\enspace}\ignorespaces}
\def\item{\par\hang\textindent} \def\op{\oplus}

\def\LH{{\bf LH}}\def\LM{{\bf LM}}\def\LT{{\bf
LT}}\def\KX{K\langle X\rangle}
\def\B{{\cal B}} \def\LC{{\bf LC}} \def\G{{\cal G}} \def\FRAC#1#2{\displaystyle{\frac{#1}{#2}}}
\def\SUM^#1_#2{\displaystyle{\sum^{#1}_{#2}}} \def\O{{\cal O}}  \def\J{{\bf J}}
\def\KS{K\langle X\rangle}\def\KXT{K\langle X,T\rangle}
{\parindent=0pt\vskip 1truecm

Let $K$ be a field. In the computational Gr\"obner basis theory for
a commutative polynomial algebra $K[x_1,...,x_n]$ or for a
noncommutative free algebra $K\langle X_1,...,X_n\rangle$, it is a
well-known fact that a homogenous Gr\"obner basis is easier to be
obtained, that is, by virtue of both the structural advantage
(mainly the degree-truncated structure) and the computational
advantage (mainly the use of a degree-preserving fast ordering),
most of the practically used commutative or noncommutative Gr\"obner
basis algorithms produce a Gr\"obner basis  by homogenizing
generators first (if the given generators are not homogeneous), and
then, in turn, producing a homogeneous Gr\"obner basis and
dehomogenizing it. On the other hand, Gr\"obner bases and the
(de)homogenization (i.e. homogenization and dehomogenization)
technique have been combined in ([11], [6], [7], [8]) to study some
global structure properties of algebras defined by relations.
Following the idea of ([11], [13], [6], [7]) concerning
(de)homogenized Gr\"obner bases, in this paper we systematize and
deepen the study on this topic. More precisely, after giving a quick
introduction to Gr\"obner bases for ideals in an algebra with a skew
multiplicative $K$-basis in Section 1, we employ the
(de)homogenization technique as used in loc. cit. to clarify through
Section 2 and Section 3 {\parindent=.5truecm\par
\item{$\bullet$} the relation between Gr\"obner bases in $R$ and homogeneous
Gr\"obner bases in $R[t]$ (Theorem 2.3, Theorem 2.5), where
$R=\oplus_{p\in\NZ}R_p$ is an $\NZ$-graded $K$-algebra with an SM
$K$-basis (i.e. a skew multiplicative $K$-basis, see the definition
in Section 1) consisting of homogeneous elements such that $R$ has a
Gr\"obner basis theory, and $R[t]$ is the polynomial ring in
commuting variable $t$ over $R$; \par
\item{$\bullet$} and the relation between Gr\"obner bases in $\KS$ and homogeneous
Gr\"obner bases in $\KXT$ (Theorem 3.3, Theorem 3.5), where $\KX
=K\langle X_1,...,X_n\rangle$ is the free $K$-algebra of $n$
generators and $\KXT =K\langle X_1,...,X_n,T\rangle$ is the fee
$K$-algebra of $n+1$ generators.\parindent=0pt\par This makes a
solid foundation for us to achieve the following goas:  Firstly, in
both cases mentioned above, demonstrate a general algorithmic
principle of obtaining a Gr\"obner basis for an ideal $I$ generated
by non-homogeneous elements and thereby obtaining a homogeneous
Gr\"obner basis for the homogenization ideal of $I$, by passing to
dealing with the homogenized generators (Proposition 2.7,
Proposition 3.7);  Secondly, find all homogeneous Gr\"obner bases in
$R[t]$ that correspond bijectively to all Gr\"obner bases in $R$
(Theorem 2.9), respectively find all homogeneous Gr\"obner bases in
$\KXT$ that correspond bijectively to all Gr\"obner bases in $\KS$
(Theorem 3.9); Thirdly, characterize all graded ideals in $R[t]$
that correspond bijectively to all ideals in $R$, respectively
characterize all graded ideals in $\KXT$ that correspond bijectively
to all ideals in $\KS$, in terms of Gr\"obner bases (Theorem 2.10,
Theorem 3.10). Based on the results obtained in previous sections,
we show in Section 4 that algebras defined by dh-closed homogeneous
Gr\"obner bases (see the definition in Section 2 and Section 3) can
be studied as Rees algebras (defined by grading filtration)
effectively via studying algebras with simpler defining relations as
demonstrated in ([7], [8]). A typical stage that may bring
Proposition 4.1 and Proposition 4.2 into play is indicated by
Theorem 2.12 (where $R[t]$ is replaced by the commutative polynomial
$K$-algebra $K[x_1,...,x_n]$ with
$R=K[x_1,...,x_{i-1},x_{i+1},...,x_n]$, $t=x_i$, $1\le i\le n$) and
Theorem 3.12 (where $\KXT =K\langle X_1,...,X_n\rangle$ with $\KX
=K\langle X_1,...,X_{i-1},X_{i+1},...,X_n\rangle$, $T=X_i$, $1\le
i\le n$).}} \v5

Throughout the paper, $\NZ$ denotes the additive monoid of nonnegative integers;
all algebras considered are associative algebras with
multiplicative identity 1; and all ideals considered are meant two-sided
ideals. If $S$ is a nonempty subset of an algebra, then we use
$\langle S\rangle$ to denote the two-sided ideal generated by $S$.
Moreover, if $K$ is a field, then we write $K^*$ for the set of
nonzero elements of $K$, i.e., $K^*=K-\{ 0\}$. \v5

\section*{1. Gr\"obner Bases w.r.t. SM $K$-bases }
Let $K$ be a field. In this section we sketch the Gr\"obner basis
theory  for $K$-algebras with an SM $K$-basis (i.e. a skew
multiplicative $K$-basis) introduced in [7]. \v5

Let $R$ be a $K$-algebra , and let $\B$ be a $K$-basis of $R$.
Adopting the commonly used notation and terminology in computational
algebra, we use lowercase letters $w,~u, ~v,~s,...$ to denote
elements in $\B$ and call an element $w\in\B$ a {\it monomial}. If
$\prec$ is a well-ordering on $\B$, $f\in R$ and
$$f=\sum^s_{i=1}\lambda_iu_i,\quad \lambda_i\in K^*,~u_i\in\B
,~u_1\prec u_2\prec\cdots\prec u_s,$$
then the {\it leading monomial}, {\it leading coefficient} and {\it
leading term} of $f$ are respectively denoted by
$$\LM (f)=u_s,\quad \LC (f)=\lambda_s,\quad \LT (f)=\lambda_su_s.$$
Furthermore, if $\prec$ satisfies the following
conditions:{\parindent=0pt\par
~(Mo1) if $w,u,v\in\B$, $u\prec v$, and $uw\ne 0$, $vw\ne 0$, then
$\LM (uw)\prec\LM (vw)$;\par
~(Mo2) if $w,u,v\in\B$, $u\prec v$, and $wu\ne 0$, $wv\ne 0$, then
$\LM (wu)\prec\LM (wv)$;\par
~(Mo3) if $w,u,v\in\B$ and $\LM (uw)=v$, then $u\prec v$, $w\prec
v$, \par
then $\prec$ is called a {\it monomial ordering on $\B$}. }\par

If $R$ has a $K$-basis $\B$ satisfying
$$u,v\in\B~\hbox{implies}~\left\{\begin{array}{l} u\cdot
v=\lambda w~\hbox{where}~\lambda \in K^*,~w\in\B,\\ \hbox{or}~u\cdot
v=0.\end{array}\right.$$ then $\B$ is called a {\it skew
multiplicative $K$-basis} (abbreviated SM $K$-basis).  If a
$K$-algebra $R$ has an SM $K$-basis and a monomial ordering $\prec$
on $\B$, then the pair $(\B ,\prec )$ is called an {\it admissible
system} of $R$; in this case the division of monomials in $R$ is
defined as follows: $u,v\in\B$, $u|v$ if and only if there is some
$\lambda\in K^*$ and $w, s\in \B$ such that $v=\lambda wus$;
furthermore, the division of monomials induces a $\prec$-compatible
division algorithm for elements in $R$, and consequently, a
Gr\"obner basis theory for $R$ may be carried out, that is, if $I$
is a nonzero ideal of $R$, then $I$ has a (finite or infinite) {\it
Gr\"obner basis} $\G$ in the sense that if $f\in I$, $f\ne 0$, then
there is some $g\in\G$ such that $\LM (g)|\LM (f)$ (see Proposition
1.2 below). \par

Commutative polynomial $K$-algebra, noncommutative free $K$-algebra,
path algebra over $K$, the coordinate algebra of a quantum affine
$n$-space over $K$, and exterior $K$-algebra are typical
$K$-algebras with an SM $K$-basis and a Gr\"obner basis theory (cf.
[1], [2], [3], [4], [5], [12]). {\parindent=0pt \v5

{\bf 1.1. Theorem} Suppose that the $K$-algebra $R$ has an
admissible system $(\B ,\prec )$ with $K$ an SM $K$-basis, and let
$I$ be an ideal of $R$. For a subset $\G\subset I$, the following
statements are equivalent:\par (i) $\G$ is a Gr\"obner basis of
$I$;\par (ii) For $f\in I$, if $f\neq 0$, then $f$ has a Gr\"obner
presentation, i.e.
$$f=\sum_{i,j}\lambda_{ij}w_{ij}g_jv_{ij},~\lambda_{ij}\in
K^*,~w_{ij},v_{ij}\in{\cal B},~g_j\in{\cal G},$$ satisfying ${\bf
LM}(w_{ij}g_jv_{ij})\preceq{\bf LM}(f)$, and there is some $j^*$
such that ${\bf LM}(w_{ij^*}g_{j^{*}}v_{ij^*})={\bf LM}(f)$;\par
(iii) $\langle{\bf LM}(I)\rangle=\langle{\bf LM}(\cal G)\rangle$,
where $\langle\LM (I)\rangle$ is the ideal of $R$ generated by the
set $\LM (I)=\{ \LM (f)~|~f\in I\}$, and $\langle\LM (\G)\rangle$ is
the ideal of $R$ generated by the set $\LM (\G )=\{ \LM (f)~|~f\in
\G\}$.\par\QED} \v5

Suppose that the $K$-algebra $R$ has an admissible system $(\B
,\prec )$ with $K$ an SM $K$-basis, and let $I$ be an ideal of $R$.
If $\G$ is a Gr\"obner basis of $I$ and any proper subset of $\G$
cannot be a Gr\"obner basis, then $\G$ is called a {\it minimal
Gr\"obner basis} of $I$. It follows easily from Theorem 1.1(ii) that
$\G$ is a minimal Gr\"obner basis of $I$ if and only if $\LM
(g_1){\not |}~\LM (g_2)$ for $g_1,g_2\in\G$ with $g_1\ne g_2$.
{\parindent=0pt\v5

{\bf 1.2. Proposition} Suppose that the $K$-algebra $R$ has an SM
$K$-basis and an admissible system $(\B ,\prec )$.  The following
two statements hold.\par
(i) Every ideal $I$ of $R$ has a minimal Gr\"obner basis:
$$\G =\{ g\in I~|~\hbox{if}~g'\in I~\hbox{and }~g'\ne
g,~\hbox{then}~\LM (g'){\not |}~\LM (g)\} .$$ \par (ii) If
$R=\oplus_{p\in\NZ}R_p$ is an $\NZ$-graded $K$-algebra and $\B$
consists of $\NZ$-homogeneous elements, then every graded ideal $I$
has a minimal homogeneous Gr\"obner basis, i.e., a minimal Gr\"obner
basis consisting of homogeneous elements (it is sufficient to
consider homogeneous elements of $I$ in (i) above). \QED}\v5

By the definition of a minimal Gr\"obner basis, it is not difficult
to see that if $\G$ is any Gr\"obner basis of the ideal $I$, then
the division algorithm enables us to produce from $\G$ a minimal
Gr\"obner basis. \v5

Let $R$ be a $K$-algebra that has an SM $K$-basis $\B$ and an
admissible system $(\B ,\prec )$. Then any nonempty subset $S\subset
R$ determines a subset of monomials
$$N(S)=\{ w\in\B~|~\LM (f){\not |}~w,~f\in S\} ,$$
which is usually called the set of {\it normal monomials} in $\B$
(mod $S$). Let $I$ be an ideal of $R$. By Theorem 1.1,  it is easy
to see that if $\G$ is a Gr\"obner basis of $I$ then $N(I)=N(\G )$;
and furthermore, each $f\in R$ has an expression of the form
$f=\sum_{i,j}\lambda_{ij}w_{ij}g_jv_{ij}+r_f$, where
$\lambda_{ij}\in K^*$, $w_{ij},v_{ij}\in\B$, $g_j\in\G$, and either
$r_f=0$ or $r_f$ has a unique linear expression of the form
$r_f=\sum_{\ell}\lambda_{\ell}w_{\ell}$ with $\lambda_{\ell}\in
K^*$, $w_{\ell}\in N(\G )$.\par

We finish this section by characterizing a Gr\"obner basis $\G$ in
terms of $N(\G )$, which, in turn, gives rise to the fundamental
decomposition of the $K$-space $R$ by $I$, respectively by $\langle
\LM(I)\rangle$. {\parindent=0pt\v5

{\bf 1.3. Theorem}  Let  $I=\langle\G\rangle$ be an ideal of $R$
generated by the subset $\G$. With notation as above, the following
statements are equivalent.\par (i) $\G$ is a Gr\"obner basis of $I$.
\par                                                   (ii) Consider
the $K$-subspace spanned by $N(\G )$, denoted $K$-span$N(\G )$.
Then the $K$-space $R$ has the decomposition
$$R=I\oplus K\hbox{-span}N(\G )=\langle\LM (I)\rangle\oplus K\hbox{-span}N(\G ).$$
(iii) The canonical image $\OV{N(\G )}$ of $N(\G )$ in $\KS
/\langle\LM (I)\rangle$ and $\KS/I$ forms a $K$-basis for $\KS
/\langle\LM (I)\rangle$ and $\KS/I$ respectively. }\v5

\section*{2. Central (De)homogenized Gr\"obner Bases}
Let $R=\oplus_{p\in\NZ}R_p$ be an $\NZ$-graded algebra over a field
$K$,  Throughout this section we fix the following assumption: $R$
{\it has an SM $K$-basis $\B$ consisting of $\NZ$-homogeneous
elements}, i.e., if $w\in\B$ then $w\in R_p$ for some $p\in\NZ$.
\par

By Section 1, if $R$ has an admissible system $(\B ,\prec)$  then
every ideal of $R$ has a Gr\"obner basis. Let $R[t]$ be the
polynomial ring in commuting variable $t$ over $R$ (i.e. $rt=tr$ for
all $r\in R$). Then, as we will see soon, with respect to the mixed
$\NZ$-gradation and a suitable monomial ordering, $R[t]$ has a
Gr\"obner basis theory, in particular, every graded ideal of $R[t]$
has a homogeneous Gr\"obner basis.  By means of the central
(de)homogenization technique as used in ([11], [6], [7]), the
present section aims to clarify in detail the relation between
Gr\"obner bases in $R$ and homogeneous Gr\"obner bases in $R[t]$.
Moreover, graded ideals in $R[t]$ that correspond bijectively to all
ideals in $R$ are characterized in terms of dh-closed Gr\"obner
bases (see the definition later).  \v5

Note that $R[t]$ has the  mixed $\NZ$-gradation, that is,
$R[t]=\op_{p\in\mathbb{Z}}R[t]_p$ is an $\NZ$-graded algebra with
the degree-$p$ homogeneous part
$$R[t]_p=\left\{\left. \displaystyle{\sum_{i+j=p}}F_it^j~\right |~F_i\in R_i,
~j\ge 0
 \right \} ,\quad p\in \mathbb{N} .$$
Considering the onto ring homomorphism $\phi$: $R[t]\rightarrow R $
defined by $\phi (t)=1$, then for each $f\in R$, there exists a
homogeneous element $F\in R[t]_p$, for some $p$, such that $\phi
(F)=f$. More precisely, if $f=f_p+f_{p-1}+ \cdots +f_{p-s}$ with
$f_p\in R_p$,  $f_{p-j}\in R_{p-j}$ and $f_p\ne 0$, then
$$f^*=f_p+tf_{p-1}+\cdots +t^sf_{p-s}$$ is a homogeneous element of
degree $p$ in $R[t]_p$ satisfying $\phi (f^*)=f$. We call the
homogeneous element $f^*$ obtained this way the {\it central
homogenization} of $f$ with respect to $t$ (for the reason that $t$
is in the center of $R[t]$). On the other hand, for an element $F\in
R[t]$, we write
$$F_*=\phi (F)$$ and call it the {\it central dehomogenization} of
$F$ with respect to $t$ (again for the reason that $t$ is in the
center of $R[t]$). Hence, if $I$ is an ideal $R$, then we write
$I^*=\{ f^*~|~f\in I\}$ and call the $\NZ$-graded ideal $\langle
I^*\rangle$ generated by $I^*$ the {\it central homogenization
ideal} of $I$ in $R[t]$ with respect to $t$; and if $J$ is an ideal
of $R[t]$, then since $\phi$ is a ring epimorphism, $\phi (J)$ is an
ideal of $R$, so we write $J_*$ for $\phi (J)=\{H_*=\phi (H)~|~H\in
J\}$ and call it the {\it central dehomogenization ideal} of $J$ in
$R$ with respect to $t$. Consequently, henceforth we will also use
the notation
$({\bf\emph{J}}_*)^*=\{(h_*)^*\mid{h\in{\bf\emph{J}}}\}$.
{\parindent=0pt\v5

{\bf 2.1. Lemma} With every definition and notation made above, the
following statements hold.\par (i) For $F,G\in R[t]$,
$(F+G)_*=F_*+G_*$, $(FG)_*=F_*G_*$.\par (ii) For any $f\in R$,
$(f^*)_*=f$.\par (iii) If $F\in R[t]_p$ and if $(F_*)^*\in R[t]_q$,
then $p\ge q$ and $t^r(F_*)^*=F$ with $r=p-q$.\par (iv) If $I$ is an
ideal of $R$, then each homogeneous element $F\in \langle
I^*\rangle$ is of the form $t^rf^*$ for some $r\in\NZ$ and $f\in
I$.\par (v) If $J$ is a graded ideal of $R[t]$, then for each $h\in
J_*$ there is some homogeneous element $F\in J$ such that
$F_*=h$.\vskip 6pt {\bf Proof} By the definition of central
(de)homogenization, the verification of (i) -- (v) is
straightforward. \QED}\v5

Suppose that the given $\NZ$-graded $K$-algebra
$R=\oplus_{p\in\NZ}R_p$ has an admissible system $(\B ,\prec_{gr})$
with $\prec_{gr}$ an $\NZ$-{\it graded monomial ordering} on $\B$,
i.e., the monomial ordering $\prec_{gr}$ is determined by a
well-ordering $\prec$ on $\B$ subject to the rule: for $u,v\in\B$,
$$u\prec_{gr}v~\hbox{if}~\left\{\begin{array}{l} \hbox{deg}(u)<~\hbox{deg}(v),\\
\hbox{or}\\
\hbox{deg}(u)=~\hbox{deg}(v)~\hbox{and}~u\prec
v,\end{array}\right.$$ where deg$(~)$ denotes the degree-function on
elements of $R$ (note that elements in $\B$ are homogeneous by our
assumption). Taking the $K$-basis $\B^*=\{ t^rw~|~w\in \B
,~r\in\NZ\}$ of $R[t]$ into account, then since $\B^*$ is obviously
a skew multiplicative $K$-basis for $R[t]$, the $\NZ$-graded
monomial ordering $\prec_{gr}$ on $\B$ extends to a monomial
ordering on $\B^*$, denoted $\prec_{t\hbox{-}gr}$, as follows:
$$t^{r_1}w_1\prec_{t\hbox{-}gr}t^{r_2}w_2~\hbox{if~and ~only~if~} w_1\prec_{gr}
w_2,~\hbox{or}~w_1=w_2~\hbox{and}~r_1<r_2.$$ Thus $R[t]$ holds a
Gr\"obner basis theory  with respect to the admissible system
$(\B^*,\prec_{t\hbox{-}gr})$.\v5

It follows from the definition of $\prec_{t\hbox{-}gr}$ that
$t^r\prec_{t\hbox{-}gr}w$ for all integers $r>0$ and all $w\in\B -\{
1\}$ (if $\B$ contains the identity element 1 of $R$). Hence,
although elements of $\B^*$ are homogeneous with respect to the
mixed $\NZ$-gradation of $R[t]$,  $\prec_{t\hbox{-}gr}$ {\it is not
a graded monomial ordering} on $\B^*$. Nevertheless,  as described
in the lemma below, since  $\prec_{gr}$ is an $\NZ$-graded monomial
ordering on $\B$, leading monomials with respect to both monomial
orderings behave in a compatible way under taking the central
(de)homogenization. {\parindent=0pt\v5

{\bf 2.2. Lemma} With notation given above, the following statements
hold.\par (i) If $f\in R$, then
$$\LM (f^*)=\LM (f)~\hbox{w.r.t.}~\prec_{t\hbox{-}gr}~\hbox{on}~\B^*.$$\par
(ii) If $F$ is a nonzero homogeneous element of $R[t]$, then
$$\LM (F_*)=\LM (F)_*~\hbox{w.r.t.}~\prec_{gr}~\hbox{on}~\B .$$\par
{\bf Proof} Since the central homogenization is done with respect to
the degree of elements in $R$, that is, if $f=f_p+f_{p-1}+ \cdots
+f_{p-s}$ with $f_p\in R_p$,  $f_{p-j}\in R_{p-j}$ and $f_p\ne 0$,
then $f^*=f_p+tf_{p-1}+\cdots +t^sf_{p-s}$, the equality $\LM
(f)=\LM (f_p)=\LM(f^*)$  follows immediately from the definitions of
$\prec_{gr}$ and $\prec_{t\hbox{-}gr}$.}\par  To prove (ii), let
$F\in R[t]_p$ be a nonzero homogeneous element of degree $p$, say
$$F=\lambda t^rw+\lambda_1t^{r_1}w_1+\cdots +\lambda_st^{r_s}w_s,$$
where $\lambda,\lambda_i,\in K^*$, $r,r_i\in\NZ$,  $w,w_i\in\B$,
such that $\LM (F)=t^rw$. Noticing that $\B$ consists of
$\NZ$-homogeneous elements and $R[t]$ has the mixed $\NZ$-gradation
by the previously fixed assumption, we have
$d(t^rw)=d(t^{r_i}w_i)=p$, $1\le i\le s$. Thus $w=w_i$ will imply
$r=r_i$ and thereby $t^rw=t^{r_i}w_i$. So we may assume that $w\ne
w_i$, $1\le i\le s$. Then it follows from the definition of
$\prec_{t\hbox{-}gr}$ that $w_i\prec_{gr}w$ and $r\le r_i$, $1\le
i\le s$. Therefore $\LM (F_*)=w=\LM (F)_*$, as desired.\QED\v5

The next result is a generalization of ([11] Theorem
2.3.2).{\parindent=0pt\v5

{\bf 2.3. Theorem} With notions and notations as fixed before, let
$I$ be an ideal of $R$, and $\langle I^*\rangle$ the central
homogenization ideal of $I$ in $R[t]$ with respect to $t$. For a
subset $\G\subset I$, the following two statements are
equivalent.\par (i)  $\G$ is a Gr\"obner basis for $I$ in $R$ with
respect to the admissible system $(\B ,\prec_{gr})$;\par (ii)
$\G^*=\{ g^*~|~g\in\G\}$  is a Gr\"obner basis for $\langle
I^*\rangle$ in $R[t]$ with respect to the admissible system
$(\B^*,\prec_{t\hbox{-}gr})$. \vskip 6pt {\bf Proof} In proving the
equivalence below, without specific indication we shall use  Lemma
2.2(i) wherever it is needed.
\par (i) $\Rightarrow$ (ii) First note that $\G^*\subset \langle I^*\rangle$.
We prove further that if $F\in \langle I^*\rangle$, then $\LM
(g^*)|\LM (F)$ for some $g\in\G$. Since $\langle I^*\rangle$ is a
graded ideal, we may assume, without loss of generality, that $F$ is
a homogeneous element. So, by Lemma 2.1(iv) we have $F=t^rf^*$ for
some $f\in I$. It follows from the equality $\LM (f^*)=\LM (f)$ that
$$\LM (F)=t^r\LM (f^*)=t^r\LM (f).\eqno{(1)}$$
If $\G$ is a Gr\"obner basis for $I$, then $\LM (f)=\lambda v\LM
(g)w$ for some $\lambda\in K^*$, $g\in\G$, and $v,w\in\B$. Thus, by
(1) above we get
$$\LM (F)=t^r\LM (f)=\lambda t^rv\LM (g^*)w.$$
This shows that $\LM (g^*)|\LM (F)$, as desired.\par (ii)
$\Rightarrow$ (i) Suppose $\G^*$ is a Gr\"obner basis for the
homogenization ideal $\langle I^*\rangle$ of $I$ in $R[t]$.  Let
$f\in I$. Then $\LM (f)=\LM (f^*)=\lambda v\LM (g^*)w=\lambda v\LM
(g)w$ for some $\lambda\in K^*$, $v,w\in\B$ and $g\in\G$. This shows
that $\LM (g)|\LM (f)$, i.e., $\G$ is a Gr\"obner basis for $I$ in
$R$.\QED}\v5

We call the Gr\"obner basis $\G^*$ obtained in Theorem 2.3 the {\it
central homogenization of} $\G$ in $R[t]$ with respect to $t$, or
$\G^*$ is a {\it central homogenized Gr\"obner basis} with respect
to $t$.\v5

By Theorem 1.3, Lemma 2.2  and Theorem 2.3, we have immediately the
following corollary. {\parindent=0pt\v5

{\bf 2.4. Corollary} Let $I$ be an arbitrary ideal of $R$.  With
notation as before, if $\G$ is a Gr\"obner basis of $I$ with respect
to the data $(\B ,\prec_{gr})$, then, with respect to the data
$(\B^*, \prec_{t\hbox{-}gr})$ we have
$$N(\langle I^*\rangle ) =N(\G^*)=\{ t^rw~|~w\in N(\G),~r\in\NZ\} ,$$
that is, the set $N(\langle I^*\rangle )$ of normal monomials in
$\B^*$ (mod $\langle I^*\rangle )$  is determined by the set
$N(I)=N(\G )$ of normal monomials in $\B$ (mod $I$). Hence, the
algebra $R[t]/\langle I^*\rangle =R[t]/\langle\G^*\rangle$ has the
$K$-basis
$$\OV{N(\langle I^*\rangle )}=\{ \OV{t^rw}~|~w\in
N(I),~r\in\NZ\}.$$}\QED\v5

We may also obtain a Gr\"obner basis for an ideal $I$ of $R$ by
dehomogenizing a homogeneous Gr\"obner basis of the ideal $\langle
I^*\rangle\subset R[t]$. Below we give a more general approach to
this assertion.{\parindent=0pt \v5

{\bf 2.5. Theorem} Let $J$ be a graded ideal of $R[t]$. If
$\mathscr{G}$ is a homogeneous Gr\"obner basis of $J$ with respect
to the data $(\B^*, \prec_{t\hbox{-}gr})$, then $\mathscr{G}_*=\{
G_*~|~G\in\mathscr{G}\}$ is a Gr\"obner basis for the ideal $J_*$ in
$R$ with respect to the data $(\B ,\prec_{gr})$.\vskip 6pt {\bf
Proof} If $\mathscr{G}$ is a Gr\"obner basis of $J$, then
$\mathscr{G}$ generates $J$ and hence $\mathscr{G}_*=\phi
(\mathscr{G} )$ generates $J_*=\phi (J)$. For a nonzero $f\in J_*$,
by Lemma 2.1(v), there exists a homogeneous element $H\in J$ such
that $H_*=f$. It follows from Lemma 2.2 that
$$\LM (f)=\LM (f^*)=\LM ((H_*)^*).\eqno{(1)}$$
On the other hand, there exists some $G\in\mathscr{G}$ such that
$\LM (G)|\LM (H)$, i.e., $$\LM (H)=\lambda t^{r_1}w\LM
(G)t^{r_2}v\eqno{(2)}$$ for some $\lambda\in K^*$, $r_1,r_2\in\NZ$,
$w,v\in\B$. But by Lemma 2.1(iii) we also have $t^r(H_*)^*=H$ for
some $r\in\NZ$, and hence
$$\LM (H)=\LM (t^r(H_*)^*)=t^r\LM (H_*)^*).\eqno{(3)}$$
So, $(1)+(2)+(3)$ yields
$$\begin{array}{rcl} \lambda t^{r_1+r_2}w\LM (G)v&=&\LM (H)\\
&=&t^r\LM ((H_*)^*)\\
&=&t^r\LM (f).\end{array}$$ Dehomogenizing both sides of the above
equality, by Lemma 2.1(i) and Lemma 2.2(ii) we obtain $$\lambda w\LM
(G_*)v=\lambda w\LM (G)_*v=\LM (f).$$ This shows that $\LM (G_*)|\LM
(f)$. Therefore, $\mathscr{G}_*$ is a Gr\"obner basis for
$J_*$.\QED}\v5

We call the Gr\"obner basis $\mathscr{G}_*$ obtained in Theorem 2.5
the {\it central dehomogenization of} $\mathscr{G}$ in $R$ with
respect to $t$, or $\mathscr{G}_*$ is a {\it central dehomogenized
Gr\"obner basis} with respect to $t$.{\parindent=0pt\v5

{\bf 2.6. Corollary} Let $I$ be an ideal of $R$. If $\mathscr{G}$ is
a homogeneous Gr\"obner basis of $\langle I^*\rangle$ in $R[t]$ with
respect to the data $(\B^*, \prec_{t\hbox{-}gr})$, then
$\mathscr{G}_*=\{ g_*~|~g\in\mathscr{G}\}$ is a Gr\"obner basis for
$I$ in $R$ with respect to the data $(B ,\prec_{gr})$. Moreover,  if
$I$ is generated by the subset $F$ and $F^*\subset\mathscr{G}$, then
$F\subset\mathscr{G}_*$. \vskip 6pt {\bf Proof} Put $J=\langle
I^*\rangle$. Then since $J_*=I$, it follows from Theorem 1.5 that if
$\mathscr{G}$ is a homogeneous Gr\"obner basis of $J$ then
$\mathscr{G}_*$ is a Gr\"obner basis for $I$. The second assertion
of the theorem is clear by Lemma 1.1(ii).\QED}\v5

Let $S$ be a nonempty subset of $R$ and $I=\langle S\rangle$ the
ideal generated by $S$. Then, with $S^*=\{ f^*~|~f\in S\}$, in
general $\langle S^*\rangle\subsetneq\langle I^*\rangle$ in $R[t]$
(for instance, consider $S=\{ y^3-x-y,~y^2+1\}$ in the commutative
polynomial ring $K[x,y]$ and $S^*$ in $K[x,y,t]$ with respect to
$t$). So, from both a practical and a computational viewpoint, it is
the right place to set up the procedure of getting a Gr\"obner basis
for $I$ and hence a Gr\"obner basis for $\langle I^*\rangle$ by
producing a homogeneous Gr\"obner basis of the graded ideal $\langle
S^*\rangle$.{\parindent=0pt\v5

{\bf 2.7. Proposition} Let $I=\langle S\rangle$ be the ideal of $R$
generated by a subset $S$. Suppose that Gr\"obner bases are
algorithmically computable in $R$ and hence in $R[t]$. Then a
Gr\"obner basis for $I$ and a homogeneous Gr\"obner basis for
$\langle I^*\rangle$ may be obtained by implementing the following
procedure:\par {\bf Step 1.}  Starting with the initial subset
$S^*=\{ f^*~|~f\in S\}$ consisting of homogeneous elements, compute
a homogeneous Gr\"obner basis $\mathscr{G}$ for the graded ideal
$\langle S^*\rangle$ of $R[t]$.\par {\bf Step 2.} Noticing $\langle
S^*\rangle_*=I$, use Theorem 2.5 and dehomogenize $\mathscr{G}$ with
respect to $t$ in order to obtain the Gr\"obner basis
$\mathscr{G}_*$ for $I$.\par {\bf Step 3.} Use Theorem 2.3 and
homogenize $\mathscr{G}_*$ with respect to $t$ in order to obtain
the homogeneous Gr\"obner basis $(\mathscr{G}_*)^*$ for the graded
ideal $\langle I^*\rangle$.\par\QED}\v5

Based on Theorem 2.3 and Theorem 2.5, we proceed now to find those
homogeneous Gr\"obner bases in $R[t]$ that correspond bijectively to
all Gr\"obner bases in $R$.\v5

Considering the central (de)homogenization with respect to $t$, a
homogeneous element $F\in R[t]$ is called {\it dh-closed} if
$(F_*)^*=F$; a subset $S$ of $R[t]$ consisting of dh-closed
homogeneous elements is called a {\it dh-closed homogeneous set}; if
a dh-closed homogeneous set $\mathscr{G}$ in $R[t]$ forms a
Gr\"obner basis with respect to $\prec_{t\hbox{-}gr}$, then it is
called a {\it dh-closed homogeneous Gr\"obner basis}.\par

To better understand the dh-closed property introduced above, we
characterize a dh-closed homogeneous element as follows.
{\parindent=0pt\v5

{\bf 2.8. Lemma} With notation as before, for a homogeneous element
$F\in R[t]$, with respect to $(\B ,\prec_{gr})$ and $(\B^*
,\prec_{t\hbox{-}gr})$ the following statements are equivalent:\par
(i) $F$ is dh-closed, i.e., $(F_*)^*=F$;\par (ii) $\LM (F_*)=\LM
(F)$;\par (iii) $F$ cannot be written as $F=t^rH$ with $H$ a
homogeneous element of $R[t]$ and $r\ge 1$;\par (iv) $t{\not |}~\LM
(F)$.\vskip6pt {\bf Proof} Using Lemma 2.1 combined with the
definitions of $\prec_{t\hbox{-}gr}$ and $(F_*)^*$, the verification
of (i) $\Rightarrow$ (ii) $\Rightarrow$ (iii) $\Rightarrow$ (iv)
$\Rightarrow$ (i) is straightforward.\QED\v5

{\bf 2.9. Theorem} With respect to the systems $(\B ,\prec_{gr})$
and $(\B^*,\prec_{t\hbox{-}gr})$, there is a one-to-one
correspondence between the set of all Gr\"obner bases in $R$ and the
set of all dh-closed homogeneous Gr\"obner bases in $R[t]$:
$$\begin{array}{ccc} \left\{
\hbox{Gr\"obner bases}~\G~\hbox{in}~R\right\}&\longleftrightarrow&
\left\{\begin{array}{l} \hbox{dh-closed homogeneous}\\
\hbox{Gr\"obner bases}~\mathscr{G}~\hbox{in}~R[t]\end{array}\right\}\\
\G&\longrightarrow&\G^*\\
\mathscr{G}_*&\longleftarrow&\mathscr{G}\end{array}$$  and this
correspondence also gives rise to a bijective map between the set of
all minimal Gr\"obner bases in $R$ and the set of all dh-closed
minimal homogeneous Gr\"obner bases in $R[t]$.\vskip 6pt {\bf Proof}
Bearing the definitions of homogenization and dehomogenization with
respect to $t$ in mind, by Theorem 2.3 and Theorem  2.5 it can be
verified directly that the given rule of correspondence defines a
one-to-one map. By the definition of a minimal Gr\"obner basis, the
second assertion follows from Lemma 2.2(i) and Lemma
2.8(ii).}\QED\v5

Below we characterize the graded ideal  generated by a dh-closed
homogeneous Gr\"obner basis in $R[t]$.  {\parindent=0pt\v5

{\bf 2.10. Theorem} With notation as before, let $J$ be a graded
ideal of $R[t]$  and $\mathscr{G}$ a minimal homogeneous Gr\"obner
basis of $J$. Under $(\B ,\prec_{gr})$ and
$(\B^*,\prec_{t\hbox{-}gr})$, the following statements are
equivalent:\par (i) $\mathscr{G}$ is a dh-closed homogeneous
Gr\"obner basis;\par (ii) $J$ has the property $\langle
(J_*)^*\rangle =J$;\par (iii) The $R[t]$-module $R[t]/J$ is
$t$-torsionfree, i.e., if $\OV f=f+J\in R[t]/J$ and $\OV f\ne 0$,
then $t\OV f\ne 0$, or equivalently, $tf\not\in J$;\par (iv)
$tR[t]\cap J=tJ$.\vskip 6pt {\bf Proof} (i) $\Rightarrow$ (ii) By
Theorem 2.5, $\mathscr{G}_*$ is a Gr\"obner basis for $J_*$ in $R$
with respect to $(\B ,\prec_{gr})$. Since $\mathscr{G}$ is
dh-closed, it follows from Theorem 2.3 that $\mathscr{G}
=(\mathscr{G}_*)^*$ is a Gr\"obner basis for $\langle
(J_*)^*\rangle$. This shows that $J=\langle\mathscr{G}\rangle
=\langle (J_*)^*\rangle$. \par (ii) $\Rightarrow$ (iii) Note that
$J$ is a graded ideal and $t$ is a homogeneous element in $R[t]$. It
is sufficient to prove that $t$ does not annihilate any nonzero
homogeneous element of $R[t]/J$. Thus, assuming $F\in R[t]_p$ and
$tF\in J$, then since $J=\langle (I_*)^*\rangle$, we have
$$(F_*)^*=((tF)_*)^*\in \langle (J_*)^*\rangle=J.$$
It follows from Lemma 2.1(iii) that there exists some $r\in\NZ$ such
that $F=t^r(F_*)^*\in J$, as desired.\par (iii) $\Leftrightarrow$
(iv) Obvious.\par (iii) $\Rightarrow$ (i) Note that $\mathscr{G}$ is
a homogeneous Gr\"obner basis by the assumption. For each
$g\in\mathscr{G}$, by Lemma 2.1(iii) there is some $r\in\NZ$ such
that $t^r(g_*)^*=g$. It follows that with respect to
$\prec_{t\hbox{-}gr}$ we have $\LM (g)=t^r\LM ((g_*)^*)$. Since
$R[t]/J$ is $t$-torsionfree, if $r>0$, then $(g_*)^*\in J$. Thus,
there is some $g'\in\mathscr{G}$ such that $g'\ne g$ and $\LM
(g')|\LM ((g_*)^*)$. Hence $\LM (g')|\LM (g)$, contradicting the
assumption that $\mathscr{G}$ is a minimal Gr\"obner basis.
Therefore, we must have $r=0$, i.e., $(g_*)^*=g$. This shows that
$\mathscr{G}$ is dh-closed.\QED\v5

{\bf 2.11. Corollary} With notation as before, let $J$ be a graded
ideal of $R[t]$. If, with respect to $(\B ,\prec_{gr})$ and
$(\B^*,\prec_{t\hbox{-}gr})$,  $J$ has a dh-closed minimal
homogeneous Gr\"obner basis, then every minimal homogeneous
Gr\"obner basis of $J$ is dh-closed. \vskip 6pt {\bf Proof} This
follows from the fact that each of the properties (ii) -- (iv) in
Theorem 2.10 does not depend on the choice of the generating set for
$J$. \QED}\v5

Let $J$ be a graded ideal of $R[t]$. If $J$ has the property
mentioned in Theorem 2.10(ii), i.e., $\langle (J_*)^*\rangle =J$,
then we call $J$ a {\it dh-closed graded ideal}. This definition
generalizes the notion of a $(\phi_*)^*$-closed graded ideal
introduced in ([10], CH.III). It is easy to see that there is a
one-to-one correspondence between the set of all ideals in $R$ and
the set of all dh-closed graded ideals in $R[t]$:
$$\begin{array}{ccc} \left\{
\hbox{ideals}~I~\hbox{in}~R\right\}&\longleftrightarrow&
\left\{\hbox{dh-closed graded ideals}~J~\hbox{in}~R[t]\right\}\\
I&\longrightarrow&\langle~I^*~\rangle\\
J_*&\longleftarrow&J\end{array}$$\par

By the foregoing argument, to know wether a given graded ideal $J$
of $R[t]$ is dh-closed, it is sufficient to compute a minimal
homogeneous Gr\"obner basis $\mathscr{G}$ for $J$ (if Gr\"obner
basis is computable in $R$), and then use the definition of a
dh-closed homogeneous set or Lemma 2.8 to check whether
$\mathscr{G}$ is dh-closed. This procedure may be realized in, for
instance, commutative polynomial $K$-algebras, noncommutative free
$K$-algebras, path algebras over $K$, the coordinate algebra of a
quantum affine $n$-space over $K$, and exterior $K$-algebras,
because Gr\"obner bases are computable in these algebras and their
polynomial extensions.\v5

Focusing on the commutative polynomial $K$-algebra $K[x_1,...,x_n]$
in $n$ variables, the good thing is that the foregoing results can
be applied to $K[x_1,...,x_n]$ with respect to each $x_i$, $1\le
i\le n$. To see this clearly, Let us put $x_i=t$ ,
$R=K[x_1,...,x_{i-1},x_{i+1},...,x_n]$, and $K[x_1,...,x_n]=R[t]$.
Moreover, let $(\B ,\prec_{gr})$ be any fixed admissible system for
$R$, where $\prec_{gr}$ is a graded monomial ordering on the
standard $K$-basis $\B$ of $R$ with respect to a fixed (weight)
$\NZ$-gradation. Then $R[t]$ has the mixed $\NZ$-gradation and the
corresponding admissible system $(\B^*,\prec_{t\hbox{-}gr})$, where
$\prec_{t\hbox{-}gr}$ is the monomial ordering obtained by extending
$\prec_{gr}$ on the standard $K$-basis $\B^*$ of $R[t]$. Instead of
mentioning a version of each result obtained before, we highlight
the respective version of Theorem 2.5 and Theorem 2.9 in this case
as follows. {\parindent=0pt\v5

{\bf 2.12. Theorem} With the preparation made above, the following
statements hold.\par (i) For each $x_i=t$, $1\le i\le n$, if
$\mathscr{G}$ is a homogeneous Gr\"obner basis of the graded ideal
$J$ in $R[t] =K[x_1,...,x_n]$ with respect to
$(\B^*,\prec_{t\hbox{-}gr})$, then $\mathscr{G}_*=\{
g_*~|~g\in\mathscr{G}\}$ is a Gr\"obner basis for the ideal $J_*$ in
$R=K[x_1,...,x_{i-1},x_{i+1},...,x_n]$ with respect to $(\B
,\prec_{gr})$. \par (ii) For each $x_i=t$, $1\le i\le n$, there is a
one-to-one correspondence between the set of all dh-closed
homogeneous Gr\"obner bases in $R[t]=K[x_1,...,x_n]$ and the set of
all Gr\"obner bases in $R=K[x_1,...,x_{i-1},x_{n+1},...,x_n]$, under
which dh-closed minimal Gr\"obner bases correspond to minimal
Gr\"obner bases. \QED}\v5

Geometrically, Theorem 2.12 may be viewed as a Gr\"obner basis
realization of the correspondence between algebraic sets in the
projective space $\mathbb{P}^{n-1}_K$ and algebraic sets in the
affine space $\mathbb{A}^{n-1}_K$, where $n\ge 2$. \v5

\def\KXT{K\langle X, T\rangle}
\section*{3. Noncentral (De)homogenized Gr\"obner Bases}
In this section  we clarify in detail how Gr\"obner bases in the
free $K$-algebra $\KS=K\langle X_1,...,X_n\rangle$ of $n$
generators are related to homogeneous Gr\"obner bases in the free
$K$-algebra $K\langle X, T\rangle =K\langle X_1,...,X_n,T\rangle$ of
$n+1$ generators, if the noncentral (de)homogenization with respect
to $T$ is employed. Moreover, in terms of dh-closed Gr\"obner bases
(see the definition later), we characterize graded ideals of $\KXT$
that correspond bijectively to all ideals in $\KS$.\par
For a general algorithmic Gr\"obner basis theory we refer to [12].
\v5

Let $\KS$ be equipped with a fixed {\it weight $\NZ$-gradation}, say
each $X_i$ has degree $n_i>0$, $1\le i\le n$. Assigning to $T$ the
degree 1 in $\KXT$ and using the same weight $n_i$ for each  $X_i$
as in $\KS$, we get the weight $\NZ$-gradation of $K\langle
X,T\rangle$ which extends the weight $\NZ$-gradation of $\KS$. Let
$\B$ and $\T{\B}$ denote the standard $K$-bases of $\KS$ and $\KXT$
respectively. To be convenient we use lowercase letters $w,u,v,...$
to denote monomials in $\B$ as before, but use capitals $W,U,V,...$
to denote monomials in $\T B$.\par In what follows, we fix an
admissible system $(\B, \prec_{gr})$ for $\KS$, where $\prec_{gr}$
is an $\NZ$-{\it graded lexicographic ordering} on $\B$ with respect
to the fixed weight $\NZ$-gradation of $\KS$, such that
$$X_{i_1}\prec_{gr}X_{i_2}\prec_{gr}\cdots\prec_{gr}X_{i_n}.$$ Then it is not
difficult to see that $\prec_{gr}$ can be extended to an $\NZ$-{\it
graded lexicographic ordering} $\prec_{_{T\hbox{-}gr}}$ on $\T{\B}$
with respect to the fixed weight $\NZ$-gradation  of $\KXT$, such
that
$$T\prec_{_{T\hbox{-}gr}} X_{i_1}\prec_{_{T\hbox{-}gr}}X_{i_2}\prec_{_{T\hbox{-}gr}}\cdots
\prec_{_{T\hbox{-}gr}}X_{i_n},$$ and thus we get the admissible
system $(\T{\B}, \prec_{_{T\hbox{-}gr}})$ for $\KXT$. With respect
to $\prec_{gr}$ and $ \prec_{_{T\hbox{-}gr}}$ we use $\LM (~)$ to
denote taking the leading monomial of elements in $\KS$ and $\KXT$
respectively. \v5

Consider the fixed $\NZ$-graded structures $\KS
=\oplus_{p\in\NZ}\KS_p$, $\KXT =\oplus_{p\in\NZ}\KXT_p$, and the
ring epimorphism $\psi :~\KXT~\longrightarrow~\KS$ defined by $\psi
(X_i)=X_i$ and $\psi (T)=1$. Then each $f\in \KS$ is the image of
some homogeneous element in $\KXT$. More precisely, if
$f=f_{p}+f_{p-1}+\cdots+f_{p-s}$ with $f_{p}\in \KS_p$, $f_{p-j}\in
\KS_{p-j}$ and $f_{p}\neq0$, then $$\T
f=f_{p}+Tf_{p-1}+\cdots+T^sf_{p-s}$$ is a homogeneous element of
degree $p$ in $\KXT_p$ such that $\psi (\T f)=f$. We call the
homogeneous element $\T f$ obtained this way the {\it noncentral
homogenization} of $f$ with respect to $T$ (for the reason that $T$
is not a commuting variable). On the other hand, for $F\in \KXT$, we
write $$F_{\sim}=\psi(F)$$ and call $F_{\sim}$ the {\it noncentral
dehomogenization} of $F$ with respect to $T$ (again for the reason
that $T$ is not a commuting variable). Furthermore, if $I=\langle
S\rangle$ is the ideal of $\KS$ generated by a subset $S$, then we
define
$$ \begin{array}{l}
\T S=\{ \T f~|~f\in S\}\cup\{X_iT-TX_i~|~1\le i\le n\},\\  \T I=\{
\T f~|~f\in I\}\cup\{X_iT-TX_i~|~1\le i\le n\},\end{array}$$ and
call the graded ideal $\langle~\T I~\rangle$ generated by $\T I$ the
{\it noncentral homogenization ideal} of $I$ in $\KXT$ with respect
to $T$; while if $J$ is an ideal of $\KXT$, then since $\psi$ is a
surjective ring homomorphism, $\psi (J)$ is an ideal of $\KS$, so we
write $J_{\sim}$ for $\psi (J)=\{H_{\sim}=\psi (H)~|~H\in J\}$ and
call it the {\it noncentral dehomogenization ideal} of $J$ in $\KS$
with respect to $T$. Consequently, henceforth we will also use the
notation
$$(J_{\sim})^{\sim}=\{(h_{\sim})^{\sim}~|~h\in J\}\cup\{
X_iT-TX_i~|~1\le i\le n\} .$$\par It is straightforward to check
that with resspect to the data $(\T{\B},\prec_{_{T\hbox{-}gr}})$,
the subset $\{ X_iT-TX_i~|~1\le i\le n\}$ of $\KXT$ forms a
homogeneous Gr\"obner basis with $\LM (X_iT-TX_i)=X_iT$, $1\le i\le
n$. In the later discussion we will freely use this fact without
extra indication.{\parindent=0pt\v5

{\bf 3.1. Lemma} With notation as fixed before, the following
properties hold.\par (i) If $F,G\in \KXT$, then
$(F+G)_{\sim}=F_{\sim}+G_{\sim}$,
$(FG)_{\sim}=F_{\sim}G_{\sim}$.\par (ii) For each nonzero $f\in
\KS$, $(\T f)_{\sim}=f$.\par (iii) Let $\mathscr{C}$ be the graded
ideal of $\KXT$ generated by $\{ X_iT-TX_i~|~1\le i\le n\}$. If
$F\in \KXT_p$, then there exists an $L\in \mathscr{C}$ and a unique
homogeneous element of the form
$H=\sum_i\lambda_iT^{r_i}w_i\in\KXT_p$, where $\lambda_i\in K^*$,
$w_i\in\B$, such that $F=L+H$; moreover there is some $r\in\NZ$ such
that $T^r(H_{\sim})^{\sim}=H$, and hence $F=L+T^r(F_{\sim})^{\sim}$.
\par
(iv) Let $\mathscr{C}$ be as in (iii) above.  If $I$ is an ideal of
$\KS$, $F\in\langle ~\T I~\rangle$ is a homogeneous element, then
there exist some $L\in \mathscr{C}$, $f\in I$ and $r\in\NZ$ such
that $F=L+T^r\T f$.\par (v) If $J$ is a graded ideal of $\KXT$ and
$\{ X_iT-TX_i~|~1\le i\le n\}\subset J$, then for each nonzero $h\in
J_{\sim}$, there exists a homogeneous element
$H=\sum_i\lambda_iT^{r_i}w_i\in J$, where $\lambda_i\in K^*$,
$r_i\in\NZ$, and $w_i\in\B$, such that for some $r\in\NZ$,
$T^r(H_{\sim})^{\sim}=H$ and $H_{\sim}=h$. \vskip 6pt {\bf Proof}
(i) and (ii) follow from the definitions of noncentral
homogenization and  noncentral dehomogenization directly.\par (iii)
Since the subset $\{ X_iT-TX_i~|~1\le i\le n\}$ is a Gr\"obner basis
in  $\KXT$ with respect to $(\T{\B},\prec_{_{T\hbox{-}gr}})$, such
that  $\LM (X_iT-TX_i)=X_iT$, $1\le i\le n$,  if $F\in\KXT_p$, then
the division of $F$ by this subset yields $F=L+H$, where $L\in
\mathscr{C}$, and $H=\sum_i\lambda_iT^{r_i}w_i$ is the unique
remainder with $\lambda_i\in K^*$, $w_i\in\B$, in which each
monomial $T^rw_i$ is of degree $p$. By the definition of
$\prec_{_{T\hbox{-}gr}}$, the definitions of noncentral
homogenization and the definition of noncentral dehomogenization, it
is not difficult to see that $H$ has the desired property.\par (iv)
By (iii), $F=L+T^r(F_{\sim})^{\sim}$ with $L\in \mathscr{C}$ and
$r\in\NZ$. Since by (ii) we have  $F_{\sim}\in \langle~\T
I~\rangle_{\sim}=I$, thus $f=F_{\sim}$ is the desired element. \par
(v) Using basic properties of homogeneous element and graded ideal
in a graded ring, this follows from the foregoing (iii). \QED} \v5

As in the case of using the central (de)homogenization, before
turning to deal with Gr\"obner bases, we are also concerned about
the behavior of leading monomials under taking the noncentral
(de)homogenization. {\parindent=0pt\v5

{\bf 3.2. Lemma}  With the assumptions and notations as fixed above,
the following statements hold.\par  (i) For each nonzero $f\in \KS$,
we have
$$\LM (f)=\LM (\T f)~\hbox{w.r.t.}~\prec_{_{T\hbox{-}gr}}~\hbox{on}~\T{\B}.$$
(ii) If $F$ is a homogeneous element in $\KXT$ such that $X_iT{\not
|}~\LM (F)$ with respect to $\prec_{_{T\hbox{-}gr}}$ for all $1\le
i\le n$, then $\LM (F)=T^rw$ for some $r\in\NZ$ and $w\in\B$, such
that $$\LM (F_{\sim})=w=\LM
(F)_{\sim}~\hbox{w.r.t.}~\prec_{gr}~\hbox{on}~\B.$$\par
{\bf Proof} Since the noncentral homogenization is done with respect
to the degree of elements in $\KS$, that is, if $f=f_p+f_{p-1}+
\cdots +f_{p-s}$ with $f_p\in \KS_p$,  $f_{p-j}\in \KS_{p-j}$ and
$f_p\ne 0$, then $\T f=f_p+Tf_{p-1}+\cdots +T^sf_{p-s}$, the
equality $\LM (f)=\LM (f_p)=\LM(\T f)$  follows immediately from the
definitions of $\prec_{gr}$ and $\prec_{t\hbox{-}gr}$.}\par  To
prove (ii), let $F\in \KXT_p$ be a nonzero homogeneous element of
degree $p$. Then by the assumption $F$ may be written as
$$F=\lambda T^rw+\lambda_1T^{r_1}X_{j_1}W_1+\lambda_2T^{r_2}X_{j_2}W_2\cdots
+\lambda_sT^{r_s}X_{j_s}W_s,$$ where $\lambda,\lambda_i,\in K^*$,
$r,r_i\in\NZ$,  $w\in\B$ and $W_i\in\T{\B}$, such that $\LM
(F)=T^rw$. Since $\B$ consists of $\NZ$-homogeneous elements and the
$\NZ$-gradation of $\KS$ extends to give the $\NZ$-gradation of
$\KXT$, we have $d(T^rw)=d(T^{r_i}X_{j_i}W_i)=p$, $1\le i\le s$.
Also note that $T$ has degree 1. Thus $w=X_{j_i}W_i$ will imply
$r=r_i$ and thereby $T^rw=T^{r_i}X_{j_i}W_i$. So we may assume,
without loss of generality, that $w\ne X_{j_i}W_i$, $1\le i\le s$.
Then it follows from the definition of $\prec_{_{T\hbox{-}gr}}$ that
$r\le r_i$, $1\le i\le n$. Hence
$X_{j_i}W_i\prec_{_{T\hbox{-}gr}}w$, $1\le i\le s$. Therefore
$(X_{j_i}W_i)_{\sim}\prec_{gr}w$, $1\le i\le n$, and consequently
$\LM (F_{\sim})=w=\LM (F)_{\sim}$, as desired.\QED\v5

The next result strengthens ([11], Theorem 2.3.1) and ([7], Theorem
8.2), in particular, the proof of (i) $\Rightarrow$ (ii) given below
improves the proof of the same deduction given in [7].
{\parindent=0pt\v5

{\bf 3.3. Theorem} With the notions and notations as fixed above,
let $I=\langle\G\rangle$ be the ideal of $\KS$ generated by a subset
$\G$, and $\langle ~\T I~\rangle$ the noncentral homogenization
ideal of $I$ in $\KXT$ with respect to $T$. The following two
statements are equivalent.\par (i) $\G$ is a Gr\"obner basis of $I$
with respect to the admissible system $(\B ,\prec_{gr})$ of
$\KS$;\par (ii) $\T{\G}=\{ \T g~|~g\in \G\}\cup\{X_iT-TX_i~|~1\le
i\le n\}$ is a homogeneous Gr\"obner basis for $\langle~\T
I~\rangle$ with respect to the admissible system
$(\T{\B},\prec_{_{T\hbox{-}gr}})$ of $\KXT$.\vskip 6pt
{\bf Proof} In proving the equivalence below, without specific
indication we shall use Lemma 3.2(i) wherever it is needed.\par (i)
$\Rightarrow$ (ii) Suppose that $\G$ is a Gr\"obner basis for $I$
with respect to the data $(\B ,\prec_{gr})$. Let $F\in\langle ~\T
I~\rangle$. Then since $\langle~\T I~\rangle$  is a graded ideal,
we may assume, without loss of generality,  that $F$ is a nonzero
homogeneous element. We want to show that there is some $D\in\T{\G}$
such that $\LM (D)|\LM (F)$, and hence $\T{\G}$ is a Gr\"obner
basis.}\par Note that $\{ X_iT-TX_i~|~1\le i\le n\}\subset\T{\G}$
with $\LM (X_iT-TX_i)=X_iT$. If $X_iT|\LM (F)$ for some $X_iT$, then
we are done. Otherwise, $X_iT{\not |}~\LM (F)$ for all $1\le i\le
n$. Thus, by Lemma 3.2(ii), $\LM (F)=T^rw$ for some $r\in\NZ$ and
$w\in\B$ such that
$$\LM (F_{\sim})=w=\LM (F)_{\sim}.\eqno{(1)}$$
On the other hand, by Lemma 3.1(iv) we have $F=L+T^{q}\T f$, where
$L$ is an element in the ideal $\mathscr{C}$ generated by $\{
X_iT-TX_i~|~1\le i\le n\}$ in $\KXT$, $q\in\NZ$, and $f\in I$. It
turns out that
$$F_{\sim}=(\T f)_{\sim}=f~\hbox{and hence}~\LM (F_{\sim})=\LM (f).\eqno{(2)}$$
But since $\G$ is a Gr\"obner basis for $I$, there is some $g\in\G$
such that $\LM (g)|\LM (f)$, i.e., there are $u,v\in\B$ such that
$$\LM (f)=u\LM (g)v=u\LM (\T g)v.\eqno{(3)}$$
Combining (1), (2), and (3) above, we have
$$w=\LM (F_{\sim})=\LM (f)=u\LM (\T g)v.$$
Therefore, $\LM (\T g)|T^rw$, i.e., $\LM (\T g)|\LM (F)$, as
desired.{\parindent=0pt\par (ii) $\Rightarrow$ (i) Suppose that
$\T{\G}$ is a Gr\"obner basis of the graded ideal $\langle~\T
I~\rangle$ in $\KXT$. If $f\in I$, then since $\T f\in \T I$, there
is some $D\in\T{\G}$ such that $\LM (D)|\LM (\T f)$. Note that $\LM
(\T f)=\LM (f)$ and thus $T{\not |}~\LM (\T f)$. Hence $D=\T g$ for
some $g\in\G$, and there are $w,v\in\B$ such that
$$\LM (f)=\LM (\T f)=w\LM (\T g)v=w\LM (g)v,$$
i.e., $\LM (g)|\LM (f)$. This shows that $\G$ is a Gr\"obner basis
for $I$ in $R$. \QED}\v5

We call the Gr\"obner basis $\T{\G}$ obtained in Theorem 3.3 the
{\it noncentral homogenization of $\G$} in $\KXT$ with respect to
$T$, or $\T{\G}$ is a {\it noncentral homogenized Gr\"obner basis}
with respect to $T$.\v5

By Theorem 1.3, Lemma 3.1 and Theorem 3.3, the following corollary
is straightforward.{\parindent=0pt\v5

{\bf 3.4. Corollary} Let $I$ be an arbitrary ideal of $\KS$.  With
notation as before, if $\G$ is a Gr\"obner basis of $I$ with respect
to the data $(\B ,\prec_{gr})$, then, with respect to the data
$(\T{\B}, \prec_{_{T\hbox{-}gr}})$ we have
$$N(\langle ~\T I~\rangle ) =N(\T{\G})=\{ T^rw~|~w\in N(\G),~r\in\NZ\} ,$$
that is, the set $N(\langle~\T I~\rangle )$ of normal monomials in
$\T{\B}$ (mod $\langle~\T I~\rangle )$  is determined by the set
$N(I)=N(\G )$ of normal monomials in $\B$ (mod $I$). Hence, the
algebra $\KXT/\langle~\T I~\rangle =\KXT /\langle\T{\G}\rangle$ has
the $K$-basis
$$\OV{N(\langle~\T I~\rangle )}=\left\{\left. \OV{T^rw}~\right |~w\in
N(I),~r\in\NZ\right\}.$$}\QED\v5

As with the central (de)homogenization with respect to the commuting
variable $t$ in section 2, we may also obtain a Gr\"obner basis for
an ideal $I$ of $\KS$ by dehomogenizing a homogeneous Gr\"obner
basis of the ideal $\langle~\T I~\rangle\subset \KXT$. Below we give
a more general approach to this assertion that generalizes
essentially ([13], Theorem 5) in which, with $I=\langle 
f_1,...,f_s\rangle$ and $J=\langle \T f_1,...,\T 
f_s,~X_iT-TX_i,~1\le i\le n\rangle$, $\mathscr{G}$ is taken to be a 
{\it reduced Gr\"obner basis} and its proof depends on the 
reducibility of $\mathscr{G}$. {\parindent=0pt \v5

{\bf 3.5. Theorem} Let $J$ be a graded ideal of $\KXT$, and suppose
that $\{ X_iT-TX_i~|~1\le i\le n\}\subset J$. If $\mathscr{G}$ is a
homogeneous Gr\"obner basis of $J$ with respect to the data
$(\T{\B}, \prec_{_{T\hbox{-}gr}})$, then $\mathscr{G}_{\sim}=\{
G_{\sim}~|~G\in\mathscr{G}\}$ is a Gr\"obner basis for the ideal
$J_{\sim}$ in $\KS$ with respect to the data $(\B
,\prec_{gr})$.\vskip 6pt {\bf Proof} If $\mathscr{G}$ is a Gr\"obner
basis of $J$, then $\mathscr{G}$ generates $J$ and hence
$\mathscr{G}_{\sim}=\psi (\mathscr{G} )$ generates $J_{\sim}=\psi
(J)$. We show next that for each nonzero $h\in J_{\sim}$, there is
some $G_{\sim}\in \mathscr{G}_{\sim}$ such that $\LM (G_{\sim})|\LM
(h)$, and hence $\mathscr{G}_{\sim}$ is a Gr\"obner basis for
$J_{\sim}$.}\par Since $\{ X_iT-TX_i~|~1\le i\le n\}\subset J$, by
Lemma 3.1(v) there exists a homogeneous element $H\in J$ and some
$r\in\NZ$ such that $T^r(H_{\sim})^{\sim}=H$ and $H_{\sim}=h$. It
follows that
$$\LM (H)=T^r\LM ((H_{\sim})^{\sim})=T^r\LM (\T h)=T^r\LM (h).\eqno{(1)}$$
On the other hand, there is some $G\in\mathscr{G}$ such that $\LM
(G)|\LM (H)$, i.e., there are $W,V\in\T{\\B}$ such that
$$\LM (H)=W\LM (G)V.\eqno{(2)}$$
But by the above (1) we must have $\LM (G)=T^qw$ for some $q\in\NZ$
and $w\in\B$. Thus, by Lemma 3.2(ii),
$$\LM (G_{\sim})=w=\LM
(G)_{\sim}~\hbox{w.r.t.}~\prec_{gr}~\hbox{on}~\B.\eqno{(3)}$$
Combining (1), (2), and (3) above, we then obtain
$$\begin{array}{rcl} \LM (h)&=&\LM (H)_{\sim}\\
&=&(W\LM (G)V)_{\sim}\\
&=&W_{\sim}\LM (G)_{\sim}V_{\sim}\\
&=&W_{\sim}\LM (G_{\sim})V_{\sim}.\end{array}$$ This shows that $\LM
(G_{\sim})|\LM (h)$, as expected.\QED\v5

We call the Gr\"obner basis $\mathscr{G}_{\sim}$ obtained in Theorem
3.5 the {\it noncentral dehomogenization of} $\mathscr{G}$ in $\KS$
with respect to $T$, or $\mathscr{G}_{\sim}$ is a {\it noncentral
dehomogenized Gr\"obner basis} with respect to
$T$.{\parindent=0pt\v5

{\bf 3.6. Corollary} Let $I$ be an ideal of $\KS$. If $\mathscr{G}$
is a homogeneous Gr\"obner basis of $\langle~\T I~\rangle$ in $\KXT$
with respect to the data $(\T{\B}, \prec_{_{T\hbox{-}gr}})$, then
$\mathscr{G}_{\sim}=\{ g_{\sim}~|~g\in\mathscr{G}\}$ is a Gr\"obner
basis for $I$ in $\KS$ with respect to the data $(B ,\prec_{gr})$.
Moreover,  if $I$ is generated by the subset $F$ and $\T
F\subset\mathscr{G}$, then $F\subset\mathscr{G}_{\sim}$. \vskip 6pt
{\bf Proof} Put $J=\langle~\T I~\rangle$. Then since $J_{\sim}=I$,
it follows from Theorem 3.5 that if $\mathscr{G}$ is a homogeneous
Gr\"obner basis of $J$ then $\mathscr{G}_{\sim}$ is a Gr\"obner
basis for $I$. The second assertion of the theorem is clear by Lemma
3.1(ii).\QED}\v5

Let $S$ be a nonempty subset of $\KS$ and $I=\langle S\rangle$ the
ideal generated by $S$. Then, with $\T S=\{ \T f~|~f\in S\}\cup\{
X_iT-TX_i~|~1\le i\le n\}$, in general $\langle~\T
S~\rangle\subsetneq\langle~\T I~\rangle$ in $\KXT$ (for instance,
consider $S=\{ Y^3-XY-X-Y,~Y^2-X+3\}$ in the free algebra $K\langle
X,Y\rangle$ and $\T S$ in $K\langle X,Y,T\rangle$ with respect to
$T$). Again, as we did in the case dealing with (de)homogenized
Gr\"obner bases with respect to the commuting variable $t$, we take
this place to set up the procedure of getting a Gr\"obner basis for
$I$ and hence a Gr\"obner basis for $\langle ~\T I~\rangle$ by
producing a homogeneous Gr\"obner basis of the graded ideal
$\langle~\T S~\rangle$.{\parindent=0pt\v5

{\bf 3.7. Proposition} Let $I=\langle S\rangle$ be the ideal of
$\KS$ as fixed above. Suppose the ground field $K$ is computable.
Then a Gr\"obner basis for $I$ and a homogeneous Gr\"obner basis for
$\langle~\T I~\rangle$ may be obtained by implementing the following
procedure:\par {\bf Step 1.}  Starting with the initial subset $$\T
S=\{ \T f~|~f\in S\}\cup\{ X_iT-TX_i~|~1\le i\le n\} ,$$ compute a
homogeneous Gr\"obner basis $\mathscr{G}$ for the graded ideal
$\langle~\T S~\rangle$ of $\KXT$.\par {\bf Step 2.} Noticing
$\langle~\T S~\rangle_{\sim}=I$, use Theorem 3.5 and dehomogenize
$\mathscr{G}$ with respect to $T$ in order to obtain the Gr\"obner
basis $\mathscr{G}_{\sim}$ for $I$.\par {\bf Step 3.} Use Theorem
3.3 and homogenize $\mathscr{G}_{\sim}$ with respect to $T$ in order
to obtain the homogeneous Gr\"obner basis
$(\mathscr{G}_{\sim})^{\sim}$ for the graded ideal $\langle~\T
I~\rangle$.\par\QED} \v5

Based on Theorem 3.3 and Theorem 3.5, we are able to determine those
homogeneous Gr\"obner bases in $\KXT$ that correspond bijectively to
all Gr\"obner bases in $\KS$.\v5

A homogeneous element $F\in \KXT$ is called {\it dh-closed} if
$(F_{\sim})^{\sim}=F$; a subset $S$ of $\KXT$ consisting of
dh-closed homogeneous elements is called a {\it dh-closed
homogeneous set}. \par

To better understand the dh-closed property introduced above for
homogeneous elements in $\KXT$, we characterize a dh-closed
homogeneous element as follows. {\parindent=0pt\v5

{\bf 3.8. Lemma} With notation as before, for a nonzero homogeneous
element $F\in\KXT$, the following two properties are equivalent with
respect to $(\B ,\prec_{gr})$ and
$(\T{\B},\prec_{_{T\hbox{-}gr}})$:\par (i) $F$ is dh-closed;\par
(ii) $F=\sum_i\lambda_iT^{r_i}w_i$ satisfying $\LM (F_{\sim})=\LM
(F)$, where $\lambda_i\in K^*$, $r_i\in\NZ$, $w_i\in\B$.\vskip6pt
{\bf Proof} This follows easily from Lemma 3.1(iii).\QED}\v5

Let $I=\langle S\rangle$ be the ideal of $\KS$ generated by a subset
$S$. Recall that we have defined $$\T S=\{ \T f~|~f\in
S\}\cup\{X_iT-TX_i~|~1\le i\le n\}$$ If $\mathscr{H}$ is a nonempty
dh-closed homogeneous set in $\KXT$ such that the subset
$$\mathscr{G}=\mathscr{H}\cup\{ X_iT-TX_i~|~1\le i\le n\}$$ forms a Gr\"obner
basis for the graded ideal $J=\langle\mathscr{G}\rangle$ with
respect to $(\T{\B},\prec_{_{T\hbox{-}gr}})$, then we call
$\mathscr{G}$ a {\it dh-closed homogeneous Gr\"obner basis}.
{\parindent=0pt\v5

{\bf 3.9. Theorem} With respect to the systems $(\B ,\prec_{gr})$
and $(\T{\B},\prec_{_{T\hbox{-}gr}})$, there is a one-to-one
correspondence between the set of all Gr\"obner bases in $\KS$ and
the set of all dh-closed homogeneous Gr\"obner bases in $\KXT$:
$$\begin{array}{ccc} \left\{
\hbox{Gr\"obner bases}~\G~\hbox{in}~\KS\right\}&\longleftrightarrow&
\left\{\begin{array}{l} \hbox{dh-closed homogeneous}\\
\hbox{Gr\"obner bases}~\mathscr{G}~\hbox{in}~\KXT\end{array}\right\}\\
\G&\longrightarrow&\T{\G}\\
\mathscr{G}_{\sim}&\longleftarrow&\mathscr{G}\end{array}$$ and this
correspondence also gives rise to a bijective map between the set of
all minimal Gr\"obner bases in $\KS$ and the set of all dh-closed
minimal homogeneous Gr\"obner bases in $\KXT$.\vskip 6pt {\bf Proof}
By the definitions of homogenization and dehomogenization with
respect to $T$, Theorem 3.3, and Theorem  3.5, it can be verified
directly that the given rule of correspondence defines a  one-to-one
map. By the definition of a minimal Gr\"obner basis, the second
assertion follows from Lemma 3.2(i) and Lemma 3.8(ii).}\QED\v5

Below we characterize the graded ideal  generated by a dh-closed
homogeneous Gr\"obner basis in $\KXT$. To make the argument more
convenient, we let $\mathscr{C}$ denote the ideal of $\KXT$
generated by the Gr\"obner basis $\{ X_iT-TX_i~|~1\le i\le n\}$, and
let $K\hbox{-span}N(\mathscr{C})$ denote the $K$-space spanned by
the set $N(\mathscr{C})$ of normal monomials in $\T{\B}$ (mod
$\mathscr{C}$). Noticing that with respect to the monomial ordering
$\prec_{_{T\hbox{-}gr}}$ on $\T{\B}$, $\LM (X_iT-TX_i)=X_iT$ for all
$1\le i\le n$, so each element $F\in K$-span$N(\mathscr{C})$ is of
the form $F=\sum_i\lambda_iT^{r_i}w_i$ with $\lambda_i\in K^*$,
$r_i\in\NZ$, and $w_i\in\B$.{\parindent=0pt\v5

{\bf 3.10. Theorem}  With the convention made above, let
$\mathscr{H}\subset K\hbox{-span}N(\mathscr{C})$ be a subset
consisting of nonzero homogeneous elements. Suppose that the subset
$\mathscr{G}=\mathscr{H}\cup\{ X_iT-TX_i~|~1\le i\le n\}$ forms a
minimal Gr\"obner basis for the graded ideal
$J=\langle\mathscr{G}\rangle$ in $\KXT$ with respect to the data
$(\T{\B}, \prec_{_{T\hbox{-}gr}})$. Then, the following statements
are equivalent:\par (i) $\mathscr{G}$ is a dh-closed homogeneous
Gr\"obner basis, i.e., $\mathscr{H}$ is a dh-closed homogeneous
set;\par (ii) $J$ has the property $\langle(J_{\sim})^{\sim}\rangle
=J$;\par (iii) $\KXT /J$ is a $T$-torsionfree (left) $\KXT$-module,
i.e., if $\OV f=f+J\in \KXT /J$ and $\OV f\ne 0$, then $T\OV f\ne
0$, or equivalently, $Tf\not\in J$;\par (iv) $TR[t]\cap J=TJ$.\vskip
6pt
{\bf Proof} (i) $\Rightarrow$ (ii)  If $\mathscr{G}$ is dh-closed,
then by Theorem 3.5, $\mathscr{G}_{\sim}$ is a Gr\"obner basis for
the ideal $J_{\sim}$ with respect to $(\B ,\prec_{gr})$.
Furthermore, it follows from Theorem 3.3 that $\mathscr{G}$ is a
Gr\"obner basis for $\langle (J_{\sim})^{\sim}\rangle$ with respect
to $(\T{\B},\prec_{_{T\hbox{-}gr}})$. Hence,
$J=\langle\mathscr{G}\rangle =\langle (J_{\sim})^{\sim}\rangle$
.\par (ii) $\Rightarrow$ (iii) Noticing that $J$ is a graded ideal
and $T$ is a homogeneous element in $\KXT$, it is sufficient to show
that $T$ does not annihilate any nonzero homogeneous element of
$\KXT /J$. Suppose $F\in\KXT_p$ and $TF\in J$. Then since $J=\langle
(J_{\sim})^{\sim}\rangle$, we have
$$(F_{\sim})^{\sim}=((TF)_{\sim})^{\sim}\in\langle (J_{\sim})^{\sim}\rangle =J.\eqno{(1)}$$
Moreover, by Lemma 3.1(iii) there exist $L\in J$ and $r\in\NZ$ such
that
$$F=L+T^r(F_{\sim})^{\sim}.\eqno{(2)}$$
Hence $(1)+(2)$ yields $F\in J$, as desired.\par (iii)
$\Leftrightarrow$ (iv) Obvious.\par
(iii) $\Rightarrow$ (i) Let $H\in \mathscr{H}-\{ X_iT-TX_i~|~1\le
i\le n\}$. Then $H=\sum_i\lambda_iT^{r_i}w_i\in
K$-span$N(\mathscr{C})$ such that $H=T^r(H_{\sim})^{\sim}$ for some
$r\in\NZ$. If $r\ge 1$, then since $\KXT /J$ is $T$-torsionfree we
must have $(H_{\sim})^{\sim}\in J$, and $\LM (H){\not |}~\LM
((H_{\sim})^{\sim})$. Hence, there exists some $H'\in \mathscr{H}-\{
H\}$ such that $\LM (H')|\LM ((H_{\sim})^{\sim})$ and consequently
$\LM (H')|\LM (H)$, contradicting the minimality of $\mathscr{G}$.
Therefore, $r=0$, i.e., $H=(H_{\sim})^{\sim}$. This shows that
$\mathscr{H}$ is dh-closed.\QED\v5

{\bf 3.11. Corollary} With notation as before, let $J$ be a graded
ideal of $\KXT$. If, with respect to $(\B ,\prec_{gr})$ and
$(\T{\B},\prec_{_{T\hbox{-}gr}})$,  $J$ has a dh-closed minimal
homogeneous Gr\"obner basis, then every minimal homogeneous
Gr\"obner basis of $J$ is dh-closed. \vskip 6pt {\bf Proof} This
follows from the fact that each of the properties (ii) -- (iv) in
Theorem 3.10 does not depend on the choice of the generating set for
$J$. \QED}\v5

Let $J$ be a graded ideal of $\KXT$. If $J$ has the property
mentioned in Theorem 3.10(ii), i.e., $\langle
(J_{\sim})^{\sim}\rangle =J$, then we call $J$ a {\it dh-closed
graded ideal}. It is easy to see that there is a one-to-one
correspondence between the set of all ideals in $\KS$ and the set of
all dh-closed graded ideals in $\KXT$:
$$\begin{array}{ccc} \left\{
\hbox{ideals}~I~\hbox{in}~\KS\right\}&\longleftrightarrow&
\left\{\hbox{dh-closed graded ideals}~J~\hbox{in}~\KXT\right\}\\
I&\longrightarrow&\langle~\T I~\rangle\\
J_{\sim}&\longleftarrow&J\end{array}$$\par
Note that in principle Gr\"obner bases are computable in $\KXT$ if
the ground field $K$ is computable. By the foregoing argument, to
know wether a given graded ideal $J$ of $\KXT$ is dh-closed, it is
sufficient to check if $J$ contains  a minimal homogeneous Gr\"obner
basis of the form $\mathscr{G}=\mathscr{H}\cup\{ X_iT-TX_i~|~1\le
i\le n\}$ in which $\mathscr{H}\subset K$-span$N(\mathscr{C})$ is a
dh-closed homogeneous set.\v5

Finally, as in the end of Section 2, let us point out that if we
start with the free $K$-algebra $K\langle X_1,...,X_n\rangle$, then
everything we have done in this section can be done with respect to
each $X_i=T$, $1\le i\le n$, that is, just work with $\KS =K\langle
X_1,...,X_{i-1},X_{i+1},...,X_n\rangle$ and $\KXT =K\langle
X_1,...,X_n\rangle$ with $T=X_i$. Also, instead of mentioning a
version of each result obtained before, we highlight the respective
version of Theorem 3.5 and Theorem 3.9 in this case as follows.
{\parindent=0pt\v5

{\bf 3.12. Theorem}  For each $X_i=T$, $1\le i\le n$, let
$\KS=K\langle X_1,...,X_{i-1},X_{i+1},...,X_n\rangle$, $\KXT
=K\langle X_1,...,X_n\rangle$ with $T=X_i$, and fix the admissible
systems $(\B ,\prec_{gr})$, $(\T{\B},\prec_{_{T\hbox{-}gr}})$ for
$\KS$ and $\KXT$ respectively, as before. The following statements
hold.\par (i) If $J$ is a graded ideal of $\KXT$ that contains the
subset $\{ X_jT-TX_j~|~j\ne i\}$ and if $\mathscr{G}$ is a
homogeneous Gr\"obner basis of $J$ with respect to
$(\T{\B},\prec_{_{T\hbox{-}gr}})$, then $\mathscr{G}_{\sim}=\{
g_{\sim}|~g\in\mathscr{G}\}$ is a Gr\"obner basis for the ideal
$J_{\sim}$ in $\KS$ with respect to  $(\B ,\prec_{gr})$. \par (ii)
There is a one-to-one correspondence between the set of all
dh-closed homogeneous Gr\"obner bases in $\KXT$  and the set of all
Gr\"obner bases in $\KS$, under which dh-closed minimal Gr\"obner
bases correspond to minimal Gr\"obner bases. \QED}\v5

\section*{4. Algebras Defined by dh-Closed Homogeneous Gr\"obner Bases}
\def\LH{{\bf LH}}
The characterization of dh-closed graded ideals in terms of
dh-closed homogeneous Gr\"obner bases given in Section 2 and Section
3 indeed provides us with an effective way to study algebras defined
by dh-homogeneous Gr\"obner bases, that is, such algebras can be
studied as Rees algebras (defined by grading filtration) via
studying algebras with simpler defining relations as demonstrated in
([9], [7], [8]). Below we present details on this conclusion.
\par All notions and notations used in previous sections are
maintained.\v5

Let $A$ be a $K$-algebra. Recall that an $\NZ$-filtration of $A$ is
a family $FA=\{ F_pA\}_{p\in\NZ}$ with each $F_pA$ a $K$-subspace of
$A$, such that (1) $1\in F_0A$; (2) $\cup_{p\in\NZ}F_pA=A$; (3)
$F_pA\subseteq F_{p+1}A$ for all $p\in\NZ$; and $F_pAF_qA\subseteq
F_{p+q}A$. If $A$ has an $\NZ$-filtration $FA$, then $FA$ determines
two $\NZ$-graded $K$-algebras $G(A)=\oplus_{p\in\NZ}G(A)_p$ with
$G(A)_p=F_pA/F_{p-1}A$, and $\T A=\oplus_{p\in\NZ}\T A_p$ with $\T
A_p=F_pA$, where $G(A)$ is called the {\it associated graded
algebra} of $A$ and $\T A$ is called the {\it Rees algebra} of
$A$.\par Let $R=\oplus_{p\in\NZ}R_p$ be an $\NZ$-graded $K$-algebra.
Then $R$ has the $\NZ$-grading filtration $FR=\{ F_pR\}_{p\in\NZ}$
with $F_pR=\oplus_{i\le p}R_i$. If $I$ is an ideal of $R$ and
$A=R/I$, then $A$ has the $\NZ$-filtration $FA=\{ F_pA\}_{p\in \NZ}$
induced by $FR$, i.e., $F_pA=(F_pR+I)/I$. For instance, if the
commutative polynomial $K$-algebra $R=K[x_1,...,x_n]$ is equipped
with the natural $\NZ$-gradation, i.e., each $x_i$ has degree 1, or
if the noncommutative free $K$-algebra $R=K\langle
X_1,...,X_n\rangle$ is equipped with the natural $\NZ$-gradation,
i.e., each $X_i$ has degree 1, then the usually used natural
$\NZ$-filtration $FA$ on $A=R/I$ is just the filtration induced by
the $\NZ$-grading filtration $FR$ of $R$. Consider the polynomial
ring $R[t]$ and the mixed $\NZ$-gradation of $R[t]$ as described in
Section 2. By ([11], [6], [7]), or in a similar way as in loc. cit.,
it can be proved that there are graded $K$-algebra isomorphisms:
$$G(A)\cong R/\langle\LH (I)\rangle ,\quad \T A\cong R[t]/\langle I^*\rangle,\eqno{(1)}$$
where $\LH (I)=\{\LH (f)~|~f\in I\}$ with $\LH (f)$ the
$\NZ$-leading homogeneous element of $f$ as defined in [7] (i.e., if
$f=f_p+f_{p-1}+\cdots f_{p-s}$ with $f_p\in R_p-\{0\}$, $f_{p-i}\in
R_{p-i}$, then $\LH (f)=f_p$), and $I^*=\{ f^*~|~f\in I\}$ with
$f^*$ the homogenization of $f$ in $R[t]$ with respect to $t$. Now,
suppose that $R$ has a Gr\"obner basis theory with respect to some
admissible system $(\B ,\prec_{gr})$ as in Section 2, and let
$J=\langle\mathscr{G}\rangle$ be a graded ideal of $R[t]$ generated
by a dh-closed homogeneous Gr\"obner basis $\mathscr{G}$. Let $I$
denote the dehomogenization ideal $J_*$ of $J$ in $R$ with respect
to $t$, i.e. $I=J_*$. Then by Theorem 2.5 and Theorem 2.3 we have
$$I=J_*=\langle\mathscr{G}_*\rangle,\quad \langle I^*\rangle =\langle (J_*)^*\rangle
=\langle\mathscr{G}\rangle =J \eqno{(2)}$$ Furthermore, from ([11],
[6], [7]) we know that $\LH (\mathscr{G}_*)=\{\LH
(g_*)~|~g_*\in\mathscr{G}_*\}$ is a Gr\"obner basis for the graded
ideal $\langle\LH (I)\rangle$ in $R$, and so
$$\langle\LH (I)\rangle =\langle\LH (\mathscr{G}_*)\rangle .\eqno{(3)}$$
It follows from $(1)+(2)+(3)$ that we have proved the
following{\parindent=0pt\v5

{\bf 4.1. Proposition} With notation and the assumption on $R$ as
above, putting $A=R/\langle\mathscr{G}_*\rangle$, then there are
graded $K$-algebra isomorphisms:
$$G(A)\cong R/\langle\LH (\mathscr{G}_*)\rangle,\quad\T A\cong R[t]/\langle\mathscr{G}\rangle =
R[t]/J.$$ \par\QED}\v5

Thus, the algebra $R[t]/\langle\mathscr{G}\rangle =\T A$ can be
studied via studying the algebras $R/\langle\mathscr{G}_*\rangle =A$
and $R/\langle\LH (\mathscr{G}_*)\rangle =G(A)$. For instance, $\T
A$ is semiprime (prime, a domain) if and only if $A$ is semiprime
(prime, a domain); if $G(A)$ is semiprime (prime, a domain), then so
are $A$ and $\T A$; if $G(A)$ is Noetherian (artinian), then so are
$A$ and $\T A$; if $G(A)$ is Noetherian with finite global
dimension, then so are $A$ and $\T A$, etc. (see [9] for more
details). \v5

Turning to the free $K$-algebras $\KS =K\langle X_1,...,X_n\rangle$
and the free $K$-algebra $\KXT =\langle X_1,...,X_n,T\rangle$, let
the admissible system $(\B ,\prec_{gr})$ for $\KS$ and the
admissible system $(\T{\B},\prec_{_{T\hbox{-}gr}})$ for $\KXT$ be as
fixed in Section 3. {\parindent=0pt\v5

{\bf 4.2. Proposition} With the convention made above, let
$\mathscr{G}$ be a dh-closed homogeneous Gr\"obner basis in $\KXT$
with respect to the data $(\T{\B},\prec_{_{T\hbox{-}gr}})$, and put
$A=\KS /\langle\mathscr{G}_{\sim}\rangle$. Considering the
$\NZ$-filtration $FA$ of $A$ induced by the (weight) $\NZ$-grading
filtration $F\KS$ of $\KS$, then there are graded $K$-algebra
isomorphisms:
$$G(A)\cong\KS /\langle\LH (\mathscr{G}_{\sim})\rangle,\quad
\T A\cong \KXT /\langle\mathscr{G}\rangle ,$$ where $\LH
(\mathscr{G}_{\sim})=\{\LH
(g_{\sim})~|~g_{\sim}\in\mathscr{G}_{\sim}\}$ with $\LH (g_{\sim})$
the $\NZ$-leading homogeneous element of $g_{\sim}$ with respect to
the $\NZ$-gradation of $\KS$ (see an explanation above). \vskip 6pt

{\bf Proof}  To be convenient, let us put
$J=\langle\mathscr{G}\rangle$ and
$I=\langle\mathscr{G}_{\sim}\rangle$. By ([11], [6], [7]), there are
graded $K$-algebra isomorphisms:
$$G(A)\cong \KS /\langle\LH (I)\rangle ,\quad \T A\cong \KXT /\langle~\T I~\rangle .\eqno{(1)}$$
By Theorem 3.5 and Theorem 3.3 we have
$$I=\langle\mathscr{G}_{\sim}\rangle=J_{\sim},\quad \langle~\T I~\rangle
=\langle (J_{\sim})^{\sim}\rangle =\langle\mathscr{G}\rangle
=J.\eqno{(2)}$$ Furthermore, from ([11], [6], [7]) we know that $\LH
(\mathscr{G}_{\sim})$ is a Gr\"obner basis for the graded ideal
$\langle\LH (I)\rangle$ in $\KS$, and so
$$\langle\LH (I)\rangle =\langle\LH (\mathscr{G}_{\sim})\rangle .\eqno{(3)}$$
It follows from $(1)+(2)+(3)$ that the desired algebra isomorphisms
are established.\QED}\v5

Thus, as demonstrated in ([7], [8]), the algebra
$\KXT/\langle\mathscr{G}\rangle =\T A$ can be studied via studying
the algebras $\KS /\langle\mathscr{G}_{\sim}\rangle =A$, $\KS
/\langle\LH (\mathscr{G}_{\sim})\rangle =G(A)$, and the monomial
algebra $\KS /\langle\LM (\mathscr{G}_{\sim})\rangle$. The reader is
referred to loc. cit. for more details.{\parindent=0pt\v5

{\bf Remark} Note that the foregoing Theorem 2.12 and Theorem 3.12
have actually provided us with a practical stage to bring
Proposition 4.1 and Proposition 4.2 into play.}

 \v5 \centerline{References}
\parindent=.8truecm

\item{[1]} B. Buchberger, Gr\"obner bases: An algorithmic method in
polynomial ideal theory. In: {\it Multidimensional Systems Theory}
(Bose, N.K., ed.), Reidel Dordrecht, 1985, 184--232.

\re{[2]} T.~Becker and V.~Weispfenning, {\it Gr\"obner Bases},
Springer-Verlag, 1993.

\item{[3]} E. L. Green, Noncommutative Gr¡§obner bases and projective
resolutions, in: {\it Proceedings of the Euroconference
Computational Methods for Representations of Groups and Algebras},
Essen, 1997, (Michler, Schneider, eds), Progress in Mathematics,
Vol. 173, Basel, Birkha¡§user Verlag, 1999, 29--60.

\item{[4]} D. Hartley and P. Tuckey, Gr\"obner Bases in Clifford and
Grassmann Algebras, {\it J. Symb. Comput.}, 20(1995), 197--205.

\item{[5]} A.~Kandri-Rody and V.~Weispfenning, Non-commutative
Gr\"obner bases in algebras of solvable type, {\it J. Symbolic
Comput.}, 9(1990), 1--26.

\item{[6]} H. Li, {\it Noncommutative Gr\"obner Bases and
Filtered-Graded Transfer}, LNM, 1795, Springer-Verlag, 2002.

\item{[7]}  H. Li, $\Gamma$-leading homogeneous algebras and Gr\"obner
bases, {\it Advanced Lectures in Mathematics}, Vol.8, International
Press and Higher Education Press, Boston-Beijing, 2009, 155--200.

\item{[8]} H. Li, On the calculation of gl.dim$G^{\NZ}(A)$ and gl.dim$\widetilde{A}$
 by using Groebner bases, {\it Algebra Colloquium}, 16(2)(2009), 181--194.

\item{[9]} H. Li and F. Van Oystaeyen, {\it Zariskian Filtrations}, $K$-Monographs
in Mathematics, Vol. 2, Kluwer Academic Publishers, 1996.

\item{[10]} H. Li and F. Van Oystaeyen, {\it A Primer of Algebraic
Geometry}, Marcel Dekker, Inc., New York$\cdot$Basel,  2000.

\item{[11]} H. Li, Y. Wu and J. Zhang, Two applications of
noncommutative Gr\"obner bases, {\it Ann. Univ. Ferrara - Sez. VII -
Sc. Mat.}, XLV(1999), 1--24.

\item{[12]} T. Mora, An introduction to commutative and noncommutative
Gr\"obner Bases, {\it Theoretic Computer Science}, 134(1994),
131--173.

\item{[13]} P. Nordbeck, On some basic applications of Gr\"obner
bases in noncommutative polynomial rings, in: {\it Gr\"obner Bases
and Applications}, Vol. 251 of {\it LMS Lecture Note Series},
Cambridge Univ. Press, Cambridge, 1998, 463--472.

\end{document}